# Geometrical Convergence Rate for Distributed Optimization with Time-Varying Directed Graphs and Uncoordinated Step-Sizes

Qingguo Lü, Huaqing Li*

College of Electronic and Information Engineering, Southwest University, Chongqing 400715

*Corresponding author: huaqingli@hotmail.com

**Abstract:** This paper studies a class of distributed optimization algorithms by a set of agents, where each agent only has access to its own local convex objective function, and the goal of the agents is to jointly minimize the sum of all the local functions. The communications among agents are described by a sequence of time-varying directed graphs which are assumed to be uniformly strongly connected. A column stochastic mixing matrices is employed in the algorithm which exactly steers all the agents to asymptotically converge to a global and consensual optimal solution even under the assumption that the step-sizes are uncoordinated. Two fairly standard conditions for achieving the geometrical convergence rate are established under the assumption that the objective functions are strong convexity and have Lipschitz continuous gradient. The theoretical analysis shows that the distributed algorithm is capable of driving the whole network to geometrically converge to an optimal solution of the convex optimization problem as long as the uncoordinated step-sizes do not exceed some upper bound. We also give an explicit analysis for the convergence rate of our algorithm through a different approach. Finally, simulation results illustrate the feasibility of the proposed algorithm and the effectiveness of the theoretical analysis throughout this paper.

**Keywords:** Distributed optimization, multi-agent systems, time-varying directed graphs, uncoordinated step-sizes, small gain theorem

## I. INTRODUCTION

In recent years, with the rapid development of science and information technology, more and more researchers have paid their research attention to multi-agent systems and have obtained many remarkable achievements. On the one hand, multi-agent system provides a theoretical research method for modeling and analysis of complex systems. On the other hand, multi-agent system is also an important branch of distributed artificial intelligence research. As one of the most important research subjects in the field of multi-agent systems, distributed optimization problem of multi-agent systems has attracted intensive research interest over the past few years due to their wide applications in distributed formation control of multiple autonomous vehicles [1], resource allocation in peer-to-peer communication networks [2], and distributed data fusion, information processing and decision making in wireless sensor networks [3]-[5], etc. Specifically, distributed optimization framework not only avoids establishing a long-distance communication system or data fusion center, but also provides better load balance to the network. In many networked systems, multi-agent systems are applied to solve such a distributed convex optimization problem where the goal is to minimize a sum of all the local objective functions, in which each local objective function is not known or shared by other agents.

In the existing literatures, distributed optimization problems are extensively solved by the (sub)gradient descent algorithm [6]-[7], the (sub)gradient-push descent algorithm [8], the fast (sub)gradient descent algorithm [9], and the dual averaging algorithm [10]. All of the above work [6]-[10] can be regarded as the consensus-based distributed (sub)gradient descent algorithms where each agent accomplishes a consensus step and then a descent step along the local (sub)gradient direction. During this time, many valuable consensus-based (sub)gradient algorithms have been addressed. Nedic et al. [8] show that the consensus-based distributed (sub)gradient descent (DGD) algorithm converges at a rate of $O(1/\sqrt{t})$ for convex Lipshitz and possibly nonsmooth objective functions when applying a diminishing step size. This coincides with the convergence rate of the centralized subgradient descent algorithm. Then, Nedic et al. [9] propose a (sub)gradient-push descent algorithm, which drives every agent to an optimal value under a fairly standard assumption of (sub)gradient boundedness. Noting that the (sub)gradient-push descent algorithm in [9] requires no knowledge of the number of agents or the graph sequence to implement. Theoretical analysis demonstrated that the (sub)gradient-push algorithm converges at a rate of $O(\ln t/\sqrt{t})$ when using a diminishing step size. This line of work has been extended to a variety of realistic conditions for distributed optimization, such as directed [8] or random communication graph [11], stochastic (sub)gradient errors [12], heterogeneous local constraints [13]-[14], linear scaling in network size [15]. Although these algorithms are intuitive and simple, it can not be neglected that they are usually slow. Among the reasons, even if the objective functions are differentiable and strongly convex, they still need to apply a diminishing step-size to converge to a consensus solution [8]-[23]. Also, the abovementioned algorithms all require the assumption of bounded

(sub)gradient to achieve the exact optimal solution, which is another shortcoming. Furthermore, the abovementioned algorithms can be fast by using a fixed step-size, but they only converge to a neighborhood of the optimal solution set.

Recent work of Xu et al. [24] had coordinately put their sights on the analysis of a class of augmented distributed gradient method (Aug-DGM) of general linear time-invariant systems with fixed topology by applying the so-called Adapt–Then–Combine scheme. Specifically, this algorithm was proved to be convergent to an exact consensual minimizer for general convex and smooth objective functions when the uncoordinated step-size is sufficiently small. Then, Nedic et al. [25] studied the distributed optimization problems with coordinated step-size over time-varying undirected/directed graphs by combining the distributed inexact gradient method and the gradient tracking technique. The theoretical analysis showed that the distributed algorithm of Nedic et al. [25] was capable of driving the whole network to geometrically converge to an optimal solution under the assumption that the global objective function is strongly convex and has Lipschitz continuous gradient. Moreover, Nedic et al. in [26] also showed that the linear convergence rates could still be achieved in the studied algorithm for uncoordinated step-sizes, where the step-sizes do not exceed some upper bounds. The work of [26] could not be extended to directed multi-agent systems in terms of the proposed framework.

Based on the above discussions, the main focus of this paper is to establish a distributed optimization algorithm in general directed multi-agent systems, which can guarantee a consensus and geometrically converge to the optimum under uncoordinated yet bounded step-sizes. We look forward to facilitating the development of a generalized theory of distributed optimization, and our ultimate goal is to design more realistic step-sizes which are capable of adaptability and promoting the practical applications. More precisely, the main contributions of this paper can be summarized as follows: (i) A distributed algorithm is proposed to analyze the convex optimization problem of multi-agent systems over time-varying and yet uniformly strongly connected directed graphs. (ii) It is important to note that, unlike the distributed descent method or the push-sum protocol proposed in [6, 8], the proposed distributed algorithm takes account of applying uncoordinated step-sizes for local optimization and is exactly ensured to converge to the minimizer even with constant step-sizes. (iii) Theoretical analysis shows that the algorithm achieves a geometrical convergence rate as long as the uncoordinated step-sizes are smaller than an explicit upper bound and no positive lower bound is required when the objective functions are strong convexity and have Lipschitz continuous gradient. Specially, we construct linearly convergent methods along with establishment of explicit bounds on their convergence rates. (iv) Specifically, simulation results demonstrate that the algorithm has a faster convergence rates compared with the well-known distributed (sub)gradient descent (DGD) algorithm [6] and the push-sum algorithm [8].

The reminder of this paper is organized as follows. We begin in Section II where we present the notations, formulate the problem of interest, introduce the communication network, and give some useful assumptions. In Section III, we consider the optimization algorithm along with the small gain theorem and some beneficial lemmas. The convergence and convergence rate results of the proposed distributed optimization algorithm are established in Section IV. Furthermore, the effectiveness of the algorithm is testified by applying a numerical example in Section V. Finally, some conclusions and future directions are drawn in Section VI.

## II. PROBLEM FORMULATION AND ASSUMPTIONS

*A. Notations*

Some standard notations throughout the paper are shown in the following. $R, R^N$ and $R^{N \times N}$ refer to the set of real numbers, the set of $N$-dimensional real column vectors and the set of $N \times N$ real matrices, respectively. $I_N$ and $0_N$ denote the $N$-dimensional identity matrix and the $N$-dimensional zero matrix, respectively. We let $1 \in R^N$ and $0 \in R^N$ refer to the column vector with all entries being one and zero, respectively. For a matrix $W$, $W^T$ (probably a vector) and $W^{-1}$ denote its transpose and inverse, respectively. We use $<\cdot>$ to represent the inner product of two vectors. For a vector $x \in R^N$, we denote its average vector as $\bar{x} = \frac{1}{N} 1^T x$ and its consensus violation as $\breve{x} = x - \frac{1}{N} 11^T x = (I - J)x = \breve{J}x$, where $J = \frac{1}{N} 11^T$ and $\breve{J} = I - J$ are two symmetric matrices. For a vector $x \in R^N$, we use $\|x\|$ to denote the standard Euclidean norm, namely $\|x\| = \sqrt{x^T x}$ and $\|x\|_\infty$ to denote the infinite norm. We also use $\|x\|_{\breve{j}} = \sqrt{<x, \breve{J}x>}$ to denote its $\breve{J}$ weighted (semi)-norm. Since $\breve{J} = \breve{J}^T \breve{J}$, we have $\|x\|_{\breve{j}} = \|\breve{J}x\|$. We use $\|W\|$ to denote the spectral norm of matrix $W$. $\nabla f(x): R^n \to R^n$ denotes the gradient of $f(x)$. For a matrix $W$, we write $W_{ij}$ or $[W]_{ij}$ to denote its $i, j$'th entry. The notation $\text{diag}(x)$ denotes a diagonal matrix whose entry in the $i$

-th row and $i$-th column is $x_i$, and the non-diagonal elements are zeros.

*B. Problem formulation*

In this paper, we consider a network of agents labeled by $V = \{1, 2, ..., N\}$, which can only exchange information with each other via local communication. The target of the multi-agent group is to collectively solve the following convex optimization problem:

$$\min f(x), \quad f(x) = \frac{1}{N}\sum_{i=1}^{N} f_i(x), \quad x \in R^n \tag{1}$$

over a global decision vector $x$. For each agent $i$, $f_i: R^n \to R$ is the convex objective function. Assume that $f_i$ is privately known by agent $i$, and probably different. Throughout the paper, we assume that the optimization problem (1) is solvable, namely, the optimal solution set $X^* = \arg\min_{x \in R^n} f(x)$ is nonempty. We denote $f^* = f(x^*), x^* \in X^*$, the optimal value and $x^*$ an optimizer of the optimization problem (1). Specifically, the objective of this paper is to solve problem (1) in a distributed manner and require that the agents do not need to exchange the information of the objective function with each other, but only share their own states with their neighbors in each step.

*C. Communication Model*

The communication topology at time $k > 0$ among agents can be modeled as a weighted directed graph $G(k) = \{V, E(k), W(k)\}$, where $V = \{1, 2, ..., N\}$ is the set of vertices with $i$ representing $i$th vertex, $E(k) \subseteq V \times V$ is the set of edges, and $W(k) = [w_{ij}(k)] \in R^{N \times N}$ is the weighted adjacency matrix with $w_{ij}(k) \geq 0$. A directed edge denoted by a pair $(j, i)$ indicates that agent $j$ can arrive at agent $i$ or agent $i$ can receive information from agent $j$. If an edge $(j, i) \in E(k)$, then agent $j$ is called a neighbor of agent $i$ and $w_{ij}(k) > 0$. The in-neighbor set and out neighbor set of agent $i$ at time $k$ are denoted by $N_i^{in}(k) = \{j \in V \mid (j, i) \in E(k)\}$ and $N_i^{out}(k) = \{j \in V \mid (i, j) \in E(k)\}$, respectively. The in-degree and out-degree of agent $i$ can be defined by $d_i^{in}(k) = \sum_{j=1}^{N} w_{ij}(k)$ and $d_i^{in}(k) = \sum_{j=1}^{N} w_{ji}(k)$, respectively. A directed path from agent $j$ to agent $i$ is a sequence of edges $(j, i_1), (i_1, i_2), ..., (i_m, i)$ in the directed graph $G$ with different nodes $i_k, k = 1, 2, ..., m$. A directed graph is strongly connected if and only if for any two distinct agents $j$ and $i$ in the set $V$, there exists a directed path from agent $j$ to agent $i$. We here adopt the following two standard assumptions on the network communication graphs, which are standard in the analysis of distributed optimization algorithm.

**Assumption 2.1:** The time-varying directed graph sequence $G(k)$ is $B_0$-strongly connected. That is, there exists an positive integer $B_0$ such that for any $k \geq 0$, the directed graph $G(k)$ with vertices $V$ and edge set $E_{B_0}(k) = \bigcup_{s=kB_0}^{(k+1)B_0 - 1} E(s)$ is strongly connected.

**Assumption 2.2:** For any $k = 0, 1, ...$, the mixing matrix $A(k) = [A_{ij}(k)] \in R^{N \times N}$ is defined by

$$A_{ij}(k) = \begin{cases} \dfrac{1}{d_j^{\text{out}}(k) + 1} & \text{whenever } j \in N_i^{\text{in}}(k) \\ 0 & \text{otherwise} \end{cases}$$

where $d_j^{\text{out}}(k) = |N_j^{\text{out}}(k)|$ is the out-degree of agent $j$ at time $k$. Also, we can conclude that $A$ is a column stochastic matrix, i.e., $\sum_{i=1}^{N} A_{ij}(k) = 1$ for any $j \in V$.

Assumption 2.1 ensures that all agents can repeatedly interact each other through repeated communications with neighbors in the whole graph sequence $G(k)$. Particularly, this assumption is considerably weaker than requiring each $G(k)$ be strongly connected for all $k > 0$. Moreover, to facilitate the analysis of our optimization algorithm, we will make the following assumptions.

**Assumption 2.3:** For each $i = 1, ..., N$, the function $f_i$ is differentiable and has Lipschitz continuous gradients, i.e., there exists a constant $L_i \in (0, +\infty)$ such that

$$\|\nabla f_i(x) - \nabla f_i(y)\| \leq L_i \|x - y\| \text{ for any } x, y \in R^n$$

As a consequence, $f$ is $\bar{L}$-smooth with $\bar{L} = (1/N) \sum_{i=1}^{N} L_i$. We will use $\hat{L} = \max_i \{L_i\}$ in the forthcoming analysis.

**Assumption 2.4:** For each $i = 1, ..., N$, its objective function $f_i : R^n \to R$ satisfies

$$f_i(x) \geq f_i(y) + <\nabla f_i(y), x - y> + \frac{\mu_i}{2} \|x - y\|^2 \text{ for any } x, y \in R^n$$

where $\mu_i \in [0, +\infty)$ and at least one $\mu_i$ is nonzero. Moreover, we will use $\hat{\mu} = \max_i \{\mu_i\}$ and $\bar{\mu} = (1/N) \sum_{i=1}^{N} \mu_i$.

**Remark 2.1:** The distributed optimization algorithms we introduced in this paper can be redefined as Push-DIGing algorithm, which exactly incorporates the push-sum protocol [8] into the distributed inexact gradient tracking technique.

# III. OPTIMIZATION ALGORITHM

*A. The Push-DIGing Algorithm*

Based on the above notations and assumptions, we are now ready to propose our main algorithm, namely, Push-DIGing algorithm, to solve the optimization problem (1), followed by its convergence analysis. In our algorithm, each agent $i$ maintains four variables $p_i(k), s_i(k), x_i(k)$ and $y_i(k)$ at each time instant $k = 1, 2, ...$, where $p_i(k), s_i(k), y_i(k)$ are three auxiliary variables and $x_i(k)$ is a key variable used to estimate the optimal solution, denoted by $x^*$. The initialization of our Push-DIGing algorithm chooses an arbitrary $x_i(0) = p_i(0) \in R^n, y_i(0) = \nabla f_i(x_i(0)) \in R$, and $s_i(0) = 1$ for all $i = 1, ..., N$. Then, every agent $i$ at each time instant $k$ updates according to the following rules:

$$\begin{aligned}
p_i(k+1) &= A_{ii}(k)(p_i(k) - \alpha_i y_i(k)) + \sum_{j \in N_i^{\text{in}}(k)} A_{ij}(k)(p_j(k) - \alpha_j y_j(k)) \\
s_i(k+1) &= A_{ii}(k)s_i(k) + \sum_{j \in N_i^{\text{in}}(k)} A_{ij}(k)s_j(k) \\
x_i(k+1) &= p_i(k+1) / s_i(k+1) \\
y_i(k+1) &= A_{ii}(k)(y_i(k) + \nabla f_i(x_i(k+1)) - \nabla f_i(x_i(k))) + \sum_{j \in N_i^{\text{in}}(k)} A_{ij}(k)(y_j(k) + \nabla f_j(x_j(k+1)) - \nabla f_j(x_j(k)))
\end{aligned} \quad (2)$$

where the positive scalars $\alpha_i > 0, i = 1, ..., N$ are step-sizes, and the vector $\nabla f_i(x_i(k))$ is the gradient of the agent $i$'s objective function $f_i(x)$ at $x = x_i(k)$.

Let $x(t) = [x_1(t), x_2(t), \cdots, x_N(t)]^T \in R^{Nn}$, and $\nabla F(x(t)) = [\nabla f_1(x_1(t)), \nabla f_2(x_2(t)), \cdots, \nabla f_N(x_N(t))]^T$ $\in R^{Nn}$. Then the Push-DIGing Algorithm (2) can be rewritten as the following compact matrix-vector form:

$$\begin{aligned}
p(k+1) &= A(k)(p(k) - Dy(k)) \\
s(k+1) &= A(k)s(k), \; S(k+1) = \text{diag}\{s(k+1)\} \\
x(k+1) &= [S(k+1)]^{-1} p(k+1) \\
y(k+1) &= A(k)(y(k) + \nabla F(x(k+1)) - \nabla F(x(k)))
\end{aligned} \quad (3)$$

where $D$ is a diagonal matrix and $[D]_{ii} = \alpha_i$ is the constant step-size of agent $i$.

**Remark 3.1:** Throughout the paper, we only study the case $n = 1$ since the analysis method in this paper can be easily extended to multi-dimensional cases.

Here we make a simple algebraic transformation on the Push-DIGing Algorithm (3), i.e.

$$\begin{aligned}
s(k+1) &= A(k)s(k) \\
S(k+1) &= \text{diag}\{s(k+1)\} \\
x(k+1) &= R(k)(x(k) - Dh(k)) \\
h(k+1) &= R(k)h(k) + (S(k+1))^{-1} A(k)(\nabla F(x(k+1)) - \nabla F(x(k)))
\end{aligned} \quad (4)$$

where $R(k) = (S(k+1))^{-1} A(k)(S(k))$, $h(k) = (S(k))^{-1} y(k)$. Noting that, under Assumption 2.1 and 2.2, it is clear that each matrix $S(k)$ is invertible, and denote $\|S^{-1}\|_{\max} = \sup_{k \geq 0} \|S^{-1}(k)\|$ which can be proven to be bounded. Also, we can prove that $R(k)$ is actually a row stochastic matrix (see Lemma 4 of [8]).

In what follows, we will introduce the following notation

$$A_B(k) = A(k)A(k-1)...A(k+1-B)$$

for any $k = 0, 1, 2, ...$ and $B = 0, 1, 2, ...$ with $B \leq k+1$ and the special case that $A_0(k) = I$ for any $k$ and $A_B(k) = I$ for any needed $k < 0$. A crucial property of the norm of $(I - \frac{1}{N}11^T)R_B(k)$ is given in the following lemma, which comes from the properties of push-sum algorithm and can be obtained from references [27, 28].

**Lemma 3.1:** Let Assumptions 2.1 and 2.2 hold, and let $B$ be an integer satisfying $B \geq B_0$. Then, for any $k = B-1, B, ...$ and any matrix $x$ with appropriate dimensions, if $x = R_B(k)y$, we have $\|\tilde{x}\| \leq \delta \|\tilde{y}\|$ where $R_B(k) = (S(k+1))^{-1} A_B(k)(S(k+1-B))$ and $\delta = Q_1(1 - \tau^{NB_0})^{\frac{B-1}{NB_0}} < 1, Q_1 = 2N \frac{1+\tau^{-NB_0}}{1-\tau^{NB_0}}, \tau = \frac{1}{N^{2+NB_0}}$.

**Proof:** We omit the proof of Lemma 3.1 since it is almost identical to that of Lemma 13 in [25].

*B. The small gain theorem*

In order to achieve a geometric convergence of the Push-DIGing algorithm, we first introduce a preliminary result, namely, the small gain theorem. It is a somewhat unusual version derived from [25]. Moreover, the original version of the theorem has received a widely research and been extensively employed in control theory [29]. Before stating the small gain theorem, we need to give some notations.

Denote the notation $u_i$ for the infinite sequence $u_i = (u_i(0), u_i(1), u_i(2),...)$, where $u_i(k) \in R^{Nn}, \forall i = 1, 2, ..., N$. Furthermore, for any positive integer $K$, let us define

$$\|u_i\|^{\lambda, K} = \max_{k=0,...,K} \frac{1}{\lambda^k} \|u_i(k)\| \tag{5}$$

$$\|u_i\|^{\lambda} = \sup_{k \geq 0} \frac{1}{\lambda^k} \|u_i(k)\| \tag{6}$$

where the parameter $\lambda \in (0,1)$ will act as the geometric convergence parameter in our later analysis. Based on the above definition, a sufficient condition on the boundedness of $\|u_i\|^{\lambda}$ will be given in the small gain theorem that we next stated. As a fundamental result in control systems, more details about the theorem can be found in [29].

**Theorem 3.1(The small gain theorem [25]):** Suppose that $u_1,...,u_m$ are sequences such that for all $K > 0$ and each $i = 1,...,m$ we have

$$\| u_{(i \bmod m)+1} \|^{\lambda,K} \leq \gamma_i \| u_i \|^{\lambda,K} + \omega_i \tag{7}$$

where $\omega_i$ are some constant, and the nonnegative constants (gains) $\gamma_1,...,\gamma_m$ satisfy

$$\gamma_1 \gamma_2 \cdots \gamma_m < 1 \tag{8}$$

Then, we have

$$\|u_1\|^{\lambda} \leq \frac{1}{1-\gamma_1 \gamma_2 \cdots \gamma_m}(\gamma_m \gamma_{m-1} \cdots \gamma_2 \omega_1 + \gamma_m \gamma_{m-1} \cdots \gamma_3 \omega_2 + ... + \gamma_m \omega_{m-1} + \omega_m) \tag{9}$$

**Proof:** The proof is omitted here since this is straightforward to derive from [25].

**Lemma 3.2:** For any matrix sequence $u_i$ and a positive constant $\lambda \in (0,1)$, if $\|u_i\|^{\lambda}$ is bounded, then $\|u_i\|$ converges to $0$ with the geometric rate $O(\lambda^k)$.

**Proof:** Assume that $\|u_i\|^{\lambda} \leq M$, where $M$ is some arbitrary nonnegative constant, then by the definitions we obtain $\sup_{k \geq 0} \frac{1}{\lambda^k} \|u_i(k)\| \leq M$, therefore $\frac{1}{\lambda^k} \|u_i(k)\| \leq M, \forall k$. The proof then follows immediately from $\|u_i(k)\| \leq M\lambda^k, \forall k$.

Before proceeding to the main proof idea, we need to define a few quantities which will use frequently in our analysis. Consider the following additional notations

$$q(k) = x(k) - 1x^* \quad \text{for any} \quad k = 0,1,... \tag{10}$$

$$z(k) = \nabla F(x(k)) - \nabla F(x(k-1)) \quad \text{for any} \quad k = 1,2,... \tag{11}$$

where $x^* \in R$ is the optimal solution of problem (1), and the initiation $z(0) = 0$.

Consider the small gain theorem, the geometric convergence of $\|q(k)\|$ can be achieved by applying this theorem to the following circle of arrows:

$$q \xrightarrow{4} z \xrightarrow{3} h \xrightarrow{2} \begin{cases} \tilde{x} \\ y \end{cases} \xrightarrow{1} q \tag{12}$$

**Remark 3.2:** We will use the small gain theorem based on the establishment of each arrow. Specially, we need to be aware of the prerequisite that the sequences $\{\|q\|, \|z\|, \|h\|, \|\tilde{x}\|, \|y\|, \|q\|\}$ are proved to be bounded. Thus, we can conclude that all quantities in the above circle of arrows converge at an geometric rate $O(\lambda^k)$. Furthermore, to apply the small gain theorem in the following analysis, we need to require that the product of gains $\gamma_i$ is less than one, which will be achieved by seeking out an appropriate step-size matrix $D$. Now, we are ready to present the establishment of each arrow in the above circle (12).

*C. Supporting lemmas*

Before introducing the Lemma 3.3, we make some definitions only used in this lemma, which distinguish the notation used in our distributed optimization problem, algorithm and analysis. Problem (1) is redefined as follows with different notation,

$$\min_{x \in R^n} g(x) = \frac{1}{N} \sum_{i=1}^{N} g_i(x) \tag{13}$$

where each function $g_i$ satisfies Assumption 2.1 and 2.2. Consider the following inexact gradient descent on the function $g$:

$$v_{k+1} = v_k - \theta \frac{1}{N} \sum_{i=1}^{N} \nabla g_i(u_i(k)) + e_k \tag{14}$$

where $\theta$ is the step-size and $e_k$ is an additive noise. Let $v^*$ be the global optimal solution of $g$ and define

$$r_k = \| v_k - v^* \| \quad \text{for any} \quad k = 0, 1, \ldots$$

Based on the above definitions, we next introduce Lemma 3.3.

**Lemma 3.3:** Suppose that

$$\sqrt{1 - \frac{\theta \overline{\mu} \beta}{2(\beta+1)}} \leq \lambda < 1 \quad \text{and} \quad 0 < \theta < \min\{\frac{\beta+1}{\overline{\mu}\beta}, \frac{1}{\overline{L}(1+\eta)}, \frac{3}{\overline{\mu}}\}$$

where $\beta \geq 2$ and $\eta > 0$. Let Assumptions 2.1 and 2.2 hold for every function $g_i$. For the problem (13), consider the sequences $\{r_k\}$ and $\{v_k\}$ be generated by the inexact gradient descent algorithm (14). Then, for any positive integer $K$, we have

$$|r|^{\lambda,K} \leq 2r_0 + \frac{\sqrt{3-\theta\overline{\mu}}}{\lambda \theta \overline{\mu}} \| e \|^{\lambda,K} + \frac{1}{\lambda \sqrt{N}} \left( \sqrt{\frac{\hat{L}(1+\eta)}{\eta \overline{\mu}}} + \frac{\hat{\mu}}{\overline{\mu}} \beta \right) \sum_{i=1}^{N} \| v - u_i \|^{\lambda,K} \tag{15}$$

**Proof:** We refer the reader to the paper [26] for the proof of Lemma 3.3.

In the following lemma, we start with the first demonstration of the circle (12) which is grounded on the error bound of the inexact gradient descent algorithm in Lemma 3.3.

**Lemma 3.4** ($\{\|\tilde{x}\|, \|y\|\} \xrightarrow{1} \|q\|$): Let Assumptions 2.2, 2.3 and 2.4 hold. Also, assume that the parameters $\alpha$ and $\lambda$ satisfy

$$\sqrt{1 - \frac{\alpha \overline{\mu} \beta}{2(\beta+1)}} \leq \lambda < 1 \quad \text{and} \quad 0 < \alpha < \min\{\frac{\beta+1}{\overline{\mu}\beta}, \frac{1}{\overline{L}(1+\eta)}, \frac{3}{\overline{\mu}}\}$$

where $\beta \geq 2$ and $\eta > 0$ are some adjustable parameters. Then, we have that for all $K = 0, 1, \ldots$,

$$\| q \|^{\lambda,K} \le 2\sqrt{N} \| \bar{x}(0) - x^* \| + (1+\sqrt{N})\left(1 + \frac{\sqrt{N}}{\lambda}\left(\sqrt{\frac{\hat{L}(1+\eta)}{\eta\bar{\mu}}} + \frac{\hat{\mu}}{\mu}\beta\right)\right) \| \check{x} \|^{\lambda,K} + C_{y,\alpha} \| y \|^{\lambda,K} \qquad (16)$$

where the positive constant $C_{y,\alpha} = \frac{\sqrt{3 - \alpha_{\max}\bar{\mu}}}{\lambda\bar{\mu}}(1 - k_D^{-1}), k_D = \frac{\alpha_{\max}}{\alpha_{\min}}$ ($\alpha_{\max} = \max_i\{\alpha_i\}, \alpha_{\min} = \min_i\{\alpha_i\}$) if $\alpha = \alpha_{\max}$; $C_{y,\alpha} = \frac{\sqrt{3 - \bar{\alpha}\bar{\mu}}}{\bar{\alpha}\sqrt{N}\lambda\bar{\mu}}\sqrt{\sum_{i=1}^{N}(\alpha_i - \bar{\alpha})^2}$ if $\alpha = \bar{\alpha}, \bar{\alpha} = \frac{1}{N}\sum_{i=1}^{N}\alpha_i$.

**Proof:** See the Appendix.

Next, to demonstrate the second arrows in the circle (12), we proceed by showing two lemmas. One is $\| h \| \to \| \check{x} \|$ and the other is $\| h \| \to \| y \|$.

**Lemma 3.5** ($\| h \| \to \| \check{x} \|$): Let Assumptions 2.1 and 2.2 hold, and let $\lambda$ be a positive constant in $(\sqrt[B]{\delta}, 1)$, where $B$ is the constant given in Lemma 3.1. Then, we get

$$\| \check{x} \|^{\lambda,K} \le \frac{\| D \|}{(\lambda^B - \delta)}\left(\delta + Q_1\frac{\lambda - \lambda^B}{1 - \lambda}\right)\| h \|^{\lambda,K} + \frac{\lambda^B}{(\lambda^B - \delta)}\sum_{t=1}^{B}\lambda^{-(t-1)}\| \check{x}(t-1) \| \qquad (17)$$

for all $K = 0, 1, ...$, where $Q_1$ is the constant as defined in Lemma 3.1.

**Proof:** Note that $x(k+1) = R(k)(x(k) - Dh(k))$. The results can be obtained by the same argument as that illustrated in the proof of Lemma 6 in [25], thus we can achieve (17).

**Lemma 3.6** ($\| h \| \to \| y \|$): For any positive integer $K$, $\| y \|^{\lambda,K} \le \| S \|_{\max}\| h \|^{\lambda,K}$ holds, where $\| S \|_{\max}$ is the constant defined above.

**Proof:** Considering $y(k) = S(k)h(k)$, we have $\| y(k) \| = \| S(k)h(k) \| \le \| S \|_{\max}\| h(k) \|$. The result immediately follows by applying the corresponding proprieties of Euclidean norm.

The next lemma presents the establishment of the third arrows in the circle (12).

**Lemma 3.7** ($\| z \|^{\lambda,K} \to \| h \|^{\lambda,K}$): Let Assumptions 2.1-2.4 hold, let the parameter $\delta$ be as given in Lemma 3.1, and let $\lambda$ be a positive constant in $(\sqrt[B]{\delta}, 1)$. Then, we have for all $K = 0, 1, ...$,

(i) $\| h \|^{\lambda,K} \le \| \check{h} \|^{\lambda,K} + \| h \|_J^{\lambda,K}$

(ii) $\| \check{h} \|^{\lambda,K} \le Q_1\| S^{-1} \|_{\max}\| C \|_{\max}\frac{\lambda(1 - \lambda^B)}{(\lambda^B - \delta)(1 - \lambda)}\| z \|^{\lambda,K} + \frac{\lambda^B}{(\lambda^B - \delta)}\sum_{t=1}^{B}\lambda^{-(t-1)}\| \check{h}(t-1) \|$

(iii) $\| h \|_J^{\lambda,K} \le \frac{\| JR \|_{\max}}{\lambda}\| h \|^{\lambda,K} + \| S^{-1} \|_{\max}\| C \|_{\max}\| z \|^{\lambda,K}$

Specially, suppose that $\| JR \|_{\max} < \lambda < 1$. Then, followed by above three items, we finally obtain

$$\|h\|^{\lambda,K} \leq \frac{\|S^{-1}\|_{\max}\|A\|_{\max}\left(1+Q_1\frac{\lambda(1-\lambda^B)}{(\lambda^B-\delta)(1-\lambda)}\right)}{1-\frac{\|JR\|_{\max}}{\lambda}}\|z\|^{\lambda,K} + \frac{\frac{\lambda^B}{(\lambda^B-\delta)}\sum_{t=1}^{B}\lambda^{-(t-1)}\|\breve{h}(t-1)\|}{1-\frac{\|JR\|_{\max}}{\lambda}} \quad (18)$$

**Proof:** See the Appendix.

The last arrow in the circle (12) demonstrated in the following Lemma is a simple consequence of the fact that the gradient of $f$ is $L$-Lipschitz.

**Lemma 3.8** ($\|q\|^{\lambda,K} \to \|z\|^{\lambda,K}$): Let assumption 2.3 holds. Then, we have for all $K = 0,1,...$, and any $0 < \lambda < 1$,

$$\|z\|^{\lambda,K} \leq \hat{L}(1+\frac{1}{\lambda})\|q\|^{\lambda,K}$$

**Proof:** The proof procedure can imitate that of Lemma 5 in [25], and thus it is omitted.

# IV. MAIN RESULTS

Based on the circle (12) established in the previous section, next we will demonstrate a major result about the geometrical convergence rate estimate for the Push-DIGing algorithm with uncoordinated step-sizes over a time-varying directed graph sequence.

**Theorem 4.1:** Let Assumptions 2.1-2.4, and Lemma 3.1 hold. Let $\alpha_{\max} = \max_i\{\alpha_i\}$ ($\alpha_{\min} = \min_i\{\alpha_i\}$) be the largest (smallest) positive entry element of the uncoordinated step-size matrix $D$ such that

$$\alpha_{\max} \in \left(0, \min\left\{\left[\frac{(1-\delta)}{2\hat{L}(1+\sqrt{N})(1+4\sqrt{N}\sqrt{\kappa})(\delta+Q_1(B-1))} \times \frac{(1-\delta)(1-\|JR\|_{\max})-4\sqrt{3}\kappa(1-k_D^{-1})\|S\|_{\max}\|S^{-1}\|_{\max}\|A\|_{\max}(BQ_1+(1/2^B-\delta))}{\|S^{-1}\|_{\max}\|A\|_{\max}(BQ_1+(1/2^B-\delta))}\right], \frac{1}{2\overline{L}}\right\}\right)$$

where $\kappa = \hat{L}/\overline{\mu}$, and $k_D = \alpha_{\max}/\alpha_{\min}$ is the condition number of the step-size matrix $D$. Suppose that the condition number $k_D$ is selected such that

$$k_D < 1 + \frac{(\lambda^B-\delta)(1-\|JR\|_{\max})}{4\sqrt{3}\kappa\|S\|_{\max}\|S^{-1}\|_{\max}\|A\|_{\max}(BQ_1+(\lambda^B-\delta))-(\lambda^B-\delta)(1-\|JR\|_{\max})}$$

Then, the sequence $x\{k\}$ generated by the Push-DIGing algorithm with uncoordinated step-sizes converges to $1x^*$ at a global geometric rate $O(\lambda^k)$, where $\lambda \in (0,1)$ is given by

$$\lambda = \max\left\{\sqrt[B]{\delta + \frac{(GF\alpha_{\max}-HK)+\sqrt{(GF\alpha_{\max}-HK)^2+4(C+HF)GK\alpha_{\max}}}{2(C+HF)}}, \sqrt{1-\frac{\alpha_{\max}\overline{\mu}}{3}}, \|JR\|_{\max}\right\}$$

with

$$\begin{aligned}
F &= \| S^{-1} \|_{\max} \| A \|_{\max} > 0 \\
G &= 2\hat{L}(1+\sqrt{N})(1+4\sqrt{N}\sqrt{\kappa})(\delta + Q_1(B-1)) > 0 \\
C &= 1 - \| JR \|_{\max} > 0 \\
H &= -4\sqrt{3}\kappa(1-k_D^{-1}) \| S \|_{\max} < 0 \\
K &= BQ_1 \| S^{-1} \|_{\max} \| A \|_{\max} > 0
\end{aligned} \qquad (19)$$

**Proof:** It is immediately obtained from Lemma 3.4-3.8 that

i) $\| q \|^{\lambda,K} \leq \gamma_{11} \| \breve{x} \|^{\lambda,K} + \gamma_{12} \| y \|^{\lambda,K} + \omega_1$

where $\gamma_{11} = (1+\sqrt{N})\left(1 + \frac{\sqrt{N}}{\lambda}\left(\sqrt{\frac{\hat{L}(1+\eta)}{\eta\bar{\mu}}} + \frac{\hat{\mu}}{\mu}\beta\right)\right)$ and $\gamma_{12} = \frac{\sqrt{3-\alpha_{\max}\bar{\mu}}}{\lambda\bar{\mu}}(1-k_D^{-1}), \omega_1 = 2\sqrt{N} \| \bar{x}(0) - x^* \|$;

ii) $\| \breve{x} \|^{\lambda,K} \leq \gamma_{21} \| h \|^{\lambda,K} + \omega_{21}$

where $\gamma_{21} = \frac{\| D \|}{(\lambda^B - \delta)}\left(\delta + Q_1 \frac{\lambda - \lambda^B}{1-\lambda}\right)$ and $\omega_{21} = \frac{\lambda^B}{(\lambda^B - \delta)} \sum_{t=1}^{B} \lambda^{-(t-1)} \| \breve{x}(t-1) \|$;

iii) $\| y \|^{\lambda,K} \leq \gamma_{22} \| h \|^{\lambda,K} + \omega_{22}$

where $\gamma_{22} = \| S \|_{\max}$ and $\omega_{22} = 0$;

iv) $\| h \|^{\lambda,K} \leq \gamma_3 \| z \|^{\lambda,K} + \omega_3$

where $\gamma_3 = \dfrac{\| S^{-1} \|_{\max} \| A \|_{\max} (1 + Q_1 \frac{\lambda(1-\lambda^B)}{(\lambda^B - \delta)(1-\lambda)})}{1 - \frac{\| JR \|_{\max}}{\lambda}}$ and $\omega_3 = \dfrac{\frac{\lambda^B}{(\lambda^B - \delta)} \sum_{t=1}^{B} \lambda^{-(t-1)} \| \breve{h}(t-1) \|}{1 - \frac{\| JR \|_{\max}}{\lambda}}$;

v) $\| z \|^{\lambda,K} \leq \gamma_4 \| q \|^{\lambda,K} + \omega_4$

where $\gamma_4 = \hat{L}(1+\frac{1}{\lambda})$ and $\omega_4 = 0$.

Moreover, to use the small gain theorem, we must choose the step-size $\alpha_{\max}$ such that

$$(\gamma_{11}\gamma_{21} + \gamma_{12}\gamma_{22})\gamma_3\gamma_4 < 1 \qquad (20)$$

which means that

$$\begin{bmatrix} \left((1+\sqrt{N})\left(1 + \frac{\sqrt{N}}{\lambda}\left(\sqrt{\frac{\hat{L}(1+\eta)}{\eta\bar{\mu}}} + \frac{\hat{\mu}}{\mu}\beta\right)\right)\frac{\| D \|}{(\lambda^B - \delta)}\left(\delta + Q_1\frac{\lambda - \lambda^B}{1-\lambda}\right) + \frac{\sqrt{3-\alpha_{\max}\bar{\mu}}}{\lambda\bar{\mu}}(1-k_D^{-1}) \| S \|_{\max}\right) \\ \times \frac{\| S^{-1} \|_{\max} \| A \|_{\max} (1 + Q_1 \frac{\lambda(1-\lambda^B)}{(\lambda^B - \delta)(1-\lambda)})}{1 - \frac{\| JR \|_{\max}}{\lambda}} \hat{L}(1+\frac{1}{\lambda}) \end{bmatrix} < 1 \qquad (21)$$

where $\beta \geq 2$, $\eta > 0$ and other constraint conditions on parameters that occur in Lemma 3.4, 3.5, and

3.7 are stated as follows:

$$0 < \alpha_{\max} < \min\{\frac{\beta+1}{\bar{\mu}\beta}, \frac{1}{\bar{L}(1+\eta)}, \frac{3}{\bar{\mu}}\} \tag{22}$$

$$\sqrt{1 - \frac{\alpha\bar{\mu}\beta}{2(\beta+1)}} \leq \lambda < 1 \tag{23}$$

$$\sqrt[B]{\delta} < \lambda < 1 \tag{24}$$

$$\| JR \|_{\max} < \lambda < 1 \tag{25}$$

Noting that $\beta \geq 2, \eta > 0$, it follows from (22) that

$$0 < \alpha_{\max} < \frac{1}{\bar{L}(1+\eta)} \tag{26}$$

Define two specific values for the parameters $\beta = 2\hat{L}/\hat{\mu}$ and $\eta = 1$ in Lemma 3.4 to obtain some concrete (probably loose) bound on the convergence rate. Furthermore, by using $0.5 \leq \lambda < 1$ and $(1-\lambda^B)/(1-\lambda) \leq B$, from relation (21), we obtain

$$\alpha_{\max} \leq \frac{(\lambda^B - \delta)}{2\hat{L}(1+\sqrt{N})(1+4\sqrt{N}\sqrt{\kappa})(\delta + Q_1(B-1))} \times \left( \frac{(\lambda^B - \delta)(1-\| JR \|_{\max})}{\| S^{-1} \|_{\max} \| A \|_{\max} (BQ_1 + (\lambda^B - \delta))} - 4\sqrt{3}\kappa(1-k_D^{-1}) \| S \|_{\max} \right) \tag{27}$$

For convenience, we rewrite formula (27) as

$$\alpha_{\max} \leq \frac{(\lambda^B - \delta)}{G}\left( \frac{C(\lambda^B - \delta)}{F(\lambda^B - \delta) + K} + H \right) \tag{28}$$

where $F, G, C, H, K$ are defined above. In order to guarantee (21) is non-emptiness or the right hand side of (28) is always positive, we require

$$k_D < 1 + \frac{(\lambda^B - \delta)(1-\| JR \|_{\max})}{4\sqrt{3}\kappa \| S \|_{\max} \| S^{-1} \|_{\max} \| A \|_{\max} (BQ_1 + (\lambda^B - \delta)) - (\lambda^B - \delta)(1-\| JR \|_{\max})}$$

Thus, (28) further implies that

$$\lambda \geq \sqrt[B]{\delta + \frac{(GF\alpha_{\max} - HK) + \sqrt{(GF\alpha_{\max} - HK)^2 + 4(C+HF)GK\alpha_{\max}}}{2(C+HF)}} \tag{29}$$

To work out $\alpha_{\max}$ when $0.5 \leq \lambda < 1$, we need

$$\alpha_{\max} \leq \frac{(1-\delta)}{2\hat{L}(1+\sqrt{N})(1+4\sqrt{N}\sqrt{\kappa})\left(\delta+Q_1(B-1)\right)} \times$$

$$\left(\frac{(1-\delta)(1-\|JR\|_{\max}) - 4\sqrt{3}\kappa(1-k_D^{-1})\|S\|_{\max}\|S^{-1}\|_{\max}\|A\|_{\max}(BQ_1+(\frac{1}{2^B}-\delta))}{\|S^{-1}\|_{\max}\|A\|_{\max}(BQ_1+(\frac{1}{2^B}-\delta))}\right) \quad (30)$$

Using $\beta/(\beta+1) \geq 2/3$ in (23), it yields

$$\sqrt{1-\frac{\alpha_{\max}\bar{\mu}}{3}} \leq \lambda < 1 \quad (31)$$

The desired result of algorithm (2) with uncoordinated step-sizes can be immediately obtained by gathering the multiple conditions for $\alpha_{\max}$ and $\lambda$.

**Remark 4.1:** Noting that the choices of $\beta, \eta, \alpha_{\max}$, and $\lambda$ may make the bounds tighter, but here we do not pay attention to this kind of problem. We aim to give an explicit convergence rate $\lambda$ for the Push-DIGing algorithm.

**Corollary 4.2:** Suppose that all the assumptions and the definitions $F, G, C, H, K$ stated in Theorem 4.1 hold. Let $B$ be a large enough integer constant such that

$$\delta = Q_1(1-\tau^{NB_0})^{\frac{B-1}{NB_0}} < 1, \quad \text{where} \quad Q_1 = 2N\frac{1+\tau^{-NB_0}}{1-\tau^{NB_0}}, \tau = \frac{1}{N^{2+NB_0}}$$

Then, for any maximum step-size $\alpha_{\max} \in (0, \frac{(C+HF)(1-\delta)^2 + HK(1-\delta)}{G(F+K)}]$ ($\alpha_{\max} = \max_i\{\alpha_i\}$), the sequence $x\{k\}$ generated by the Push-DIGing algorithm with uncoordinated step-sizes converges to $\mathbf{1}x^*$ at a global geometric rate $O(\lambda^k)$, where

$$\lambda = \sqrt[2B]{1-\frac{\alpha_{\max}\bar{\mu}}{3}} \quad \text{if we choose} \quad \alpha_{\max} \in (0, M], \text{ and}$$

$$\lambda = \sqrt[B]{\frac{(2\delta(C+HF)-HK) + \sqrt{(2\delta(C+HF)-HK)^2 + 4(C+HF)(HK\delta + G(F+K)\alpha_{\max} - (C+HF)\delta^2)}}{2(C+HF)}}$$

If we choose $\alpha_{\max} \in (M, \frac{(C+HF)(1-\delta)^2 + HK(1-\delta)}{G(F+K)}]$.

where $M = \frac{3(1-\lambda_{mid}^{2B})}{\bar{\mu}} = \frac{(C+HF)(\lambda_{mid}^B - \delta)^2 + HK(\lambda_{mid}^B - \delta)}{G(F+K)}$ ($\lambda_{mid}$ is an intermediate variable detailedly described in the proof).

**Proof:** It directly follows from (28) that

$$\alpha_{\max} \leq \frac{(C+HF)(\lambda^B - \delta)^2 + HK(\lambda^B - \delta)}{GF(\lambda^B - \delta) + GK} \quad (32)$$

where $F, G, C, H, K$ are defined in Theorem 4.1. Recall (31) that

$$\alpha_{max} \geq \frac{3(1-\lambda^2)}{\bar{\mu}} \qquad (33)$$

Then, using relations (32) and (33), we can conclude that there exists $\lambda \in (\sqrt[B]{\delta}, 1)$ such that

$$[\frac{3(1-\lambda^2)}{\bar{\mu}}, \frac{(C+HF)(\lambda^B-\delta)^2 + HK(\lambda^B-\delta)}{GF(\lambda^B-\delta)+GK}] \neq \varnothing \qquad (34)$$

Here, we study a smaller interval by enlarging the left-side and reducing the right-side of the interval in (34). Since $B \geq 1$, we will prove that

$$[\frac{3(1-\lambda^{2B})}{\bar{\mu}}, \frac{(C+HF)(\lambda^B-\delta)^2 + HK(\lambda^B-\delta)}{G(F+K)}] \neq \varnothing \qquad (35)$$

Noting that when $\lambda$ increases from $\sqrt[B]{\delta}$ to 1, the left-side of (35) is decreasing from $\frac{3(1-\delta^2)}{\bar{\mu}}$ to 0 monotonically while the right-side is increasing from 0 to $\frac{(C+HF)(1-\delta)^2 + HK(1-\delta)}{G(F+K)}$ monotonically. Thus, when $\lambda$ varies from $\sqrt[B]{\delta}$ to 1, the critical value $M$ of the interval in (35) is valid when $\lambda$ is given by

$$\lambda = \sqrt[B]{\frac{(2\delta\bar{\mu}(C+HF) - \bar{\mu}HK) + \sqrt{(\bar{\mu}HK - 2\delta\bar{\mu}(C+HF))^2 + 4(\bar{\mu}(C+HF) + 3G(F+K))(\delta\bar{\mu}HK + 3G(F+K) - \bar{\mu}(C+HF)\delta^2)}}{2(\bar{\mu}(C+HF) + 3G(F+K))}}$$

(36)

Here, we assume that the value obtained in (36) is $\lambda_{mid}$. Thus, if we choose $\alpha_{max}$ such that $\alpha_{max} \in (0, M]$, we can use $\lambda = \sqrt[2B]{1 - \frac{\alpha_{max}\bar{\mu}}{3}}$, while for $\alpha_{max} \in \left(M, \frac{(C+HF)(1-\delta)^2 + HK(1-\delta)}{G(F+K)}\right)$, we can also set

$$\lambda = \sqrt[B]{\frac{(2\delta(C+HF) - HK) + \sqrt{(2\delta(C+HF) - HK)^2 + 4(C+HF)(HK\delta + G(F+K)\alpha_{max} - (C+HF)\delta^2)}}{2(C+HF)}}$$

The proof is thus completed.

## V. NUMERICAL SIMULATION RESULTS

In this section, a numerical example is presented to illustrate the feasibility of the proposed algorithm and the correctness of theoretical analysis throughout the paper. We consider a class of parameter estimation problem in wireless-sensor network. In the simulation, we use five sensors to cooperatively estimate a parameter $x$. Thus, the optimization problem is to find $x \in R$ to minimize $F(x) = \sum_{i=1}^{N} f_i(x) = \sum_{i=1}^{N} (a_i + \frac{\|x - c_i\|^2}{b_i})$, where $a_i$ is the initial state of local objective $f_i(x)$, $b_i$ is the control gain and $c_i$ is the observation known only to sensor $i$. Without loss of generality, we

set $a = [1, 2, 3, 4, 5], b = [3.33, 1.67, 1.11, 0.83, 0.67], c = [0.2, 0.4, 0.6, 0.8, 1]$, the uncoordinated step-sizes $D = [0.035, 0.015, 0.025, 0.045, 0.055]$, and the initial condition $x_i(0) \in (0, 1)$ for all $i = 1, \ldots, 5$.

Then, the simulation results are to be compared with the well-known distributed gradient decent algorithm (DGD) [6], Push-sum algorithm [8] in Fig. 1. It shows that the Push-DIGing algorithm we introduced in this paper has a linear convergence rates while the DGD and Push-sum algorithm only have sublinear convergence rates under the same settings.

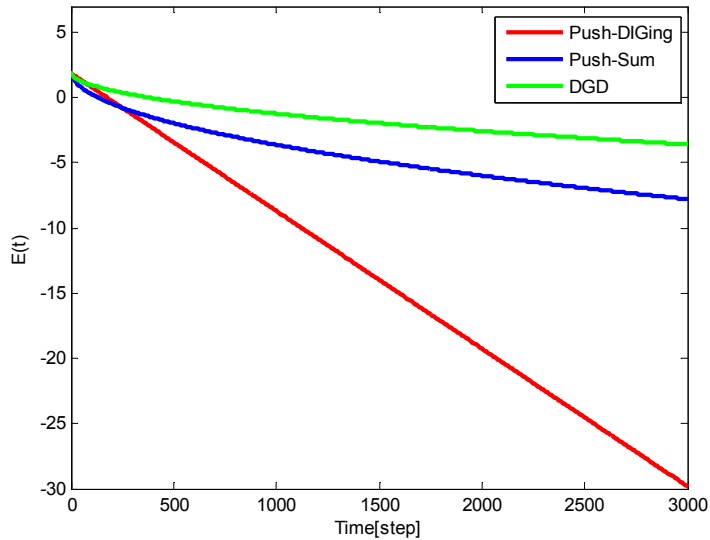

Fig. 1: Time evolution of $\log_{10} \sum_{i=1}^{5} \frac{\|x_i(k) - x^*\|}{\|x_i(0) - x^*\|}$ among Push-DIGing algorithm, Push-sum algorithm and DGD.

## VI. CONCLUSIONS AND FUTURE WORK

In this paper, the Push-DIGing algorithm for solving the convex optimization problem (1) with uncoordinated step-sizes has been studied in detail. Under some standard assumptions on the objective function and network connectivity, it has been proven that our algorithm is able to achieve a geometrical convergence rate over time-varying directed graphs. We also provided an explicit choice of $\alpha_{\max}$ and derived an explicit convergence rate $\lambda$ based on the small gain theorem. Furthermore, the correctness and effectiveness of our theoretical results are demonstrated by using a numerical example. Future work should include the case of more complex constrained convex optimization problem, and event-triggered communication among agents in the dynamic (time-varying) networks.

# APPENDIX

**Proof of Lemma 3.4:** Multiplying $\frac{1}{N}\mathbf{1}^\mathrm{T}$ on the both sides of the last equation of (3) and noting that $A(k)$ is a column stochastic matrix, we then have

$$\bar{y}(k+1) - \frac{1}{N}\mathbf{1}^\mathrm{T}\nabla F(x(k+1)) = \bar{y}(k) - \frac{1}{N}\mathbf{1}^\mathrm{T}\nabla F(x(k)) \tag{37}$$

where $\bar{y} = \frac{1}{N}\mathbf{1}^\mathrm{T}y$. Since $y_i(0) = \nabla f_i(x_i(0)), i \in V$, it follows that

$$\bar{y}(k) = \frac{1}{N}\sum_{i=1}^{N}\nabla f_i(x_i(k)) \tag{38}$$

Let us consider the evolution of $\bar{p}(k)$. Multiplying $\frac{1}{N}\mathbf{1}^\mathrm{T}$ at the both sides of the first equation of (3), one can obtain

$$\bar{p}(k+1) = \bar{p}(k) - \alpha\frac{1}{N}\sum_{i=1}^{N}\nabla f_i(x_i(k)) + (\alpha\frac{1}{N}\mathbf{1}^\mathrm{T} - \frac{1}{N}\mathbf{1}^\mathrm{T}D)y(k) \tag{39}$$

Applying Lemma 3.3 to the recursion relation of (39), and $0 < \alpha < \frac{3}{\bar{\mu}}$, we can achieve

$$\|\bar{p} - x^*\|^{\lambda,K} \leq 2\|\bar{p}(0) - x^*\| + \frac{1}{\lambda\sqrt{N}}\left(\sqrt{\frac{\hat{L}(1+\eta)}{\eta\bar{\mu}}} + \frac{\hat{\mu}}{\bar{\mu}}\beta\right)\sum_{i=1}^{N}\|\bar{p} - x_i\|^{\lambda,K}$$
$$+ \frac{\sqrt{3-\alpha\bar{\mu}}}{N\lambda\bar{\mu}}\|(\mathbf{1}^\mathrm{T} - \mathbf{1}^\mathrm{T}\frac{D}{\alpha})y\|^{\lambda,K} \tag{40}$$

Let us analyze the summation in the second term of (40). Since $x(k+1) = [S(k+1)]^{-1}s(k+1)$ and $x(0) = p(0)$, we obtain

$$\begin{aligned}\sum_{i=1}^{N}\|\bar{p} - x_i\|^{\lambda,K} &= \sum_{i=1}^{N}\|(\bar{p} - \bar{x}) + (\bar{x} - x_i)\|^{\lambda,K}\\ &\leq N\|\bar{p} - \bar{x}\|^{\lambda,K} + \sum_{i=1}^{N}\|\bar{x} - x_i\|^{\lambda,K}\\ &\leq N\|\frac{1}{N}\mathbf{1}^\mathrm{T}Sx - \frac{1}{N}\mathbf{1}^\mathrm{T}x\|^{\lambda,K} + \sum_{i=1}^{N}\|\bar{x} - x_i\|^{\lambda,K}\\ &= \|(\mathbf{1}-s)^\mathrm{T}x\|^{\lambda,K} + \sqrt{N}\|\check{x}\|^{\lambda,K}\end{aligned} \tag{41}$$

Since $q(k) = x(k) - \mathbf{1}x^* = x(k) - \mathbf{1}\bar{x}(k) + \mathbf{1}\bar{x}(k) - \mathbf{1}\bar{p}(k) + \mathbf{1}\bar{p}(k) - \mathbf{1}x^*$, it follows that

$$q(k) = \check{x}(k) + \frac{1}{N}\mathbf{1}(\mathbf{1}-s(k))^\mathrm{T}x(k) + \mathbf{1}(\bar{p}(k) - x^*) \tag{42}$$

Together with (40) and (41), (42) further implies

$$\|q\|^{\lambda,K} \leq \left(1 + \frac{\sqrt{N}}{\lambda}\sqrt{\frac{\hat{L}(1+\eta)}{\eta\bar{\mu}} + \frac{\hat{\mu}}{\mu}\beta}\right)\|\breve{x}\|^{\lambda,K} + 2\sqrt{N}\|\bar{x}(0) - x^*\|$$
$$+ \left((\sqrt{N})^{-1} + \frac{1}{\lambda}\sqrt{\frac{\hat{L}(1+\eta)}{\eta\bar{\mu}} + \frac{\hat{\mu}}{\mu}\beta}\right)\|(1-s)^T x\|^{\lambda,K} + \frac{\sqrt{3-\alpha\bar{\mu}}}{\sqrt{N}\lambda\bar{\mu}}\|(1^T - 1^T\frac{D}{\alpha})y\|^{\lambda,K}$$
(43)

Finally, let us bound the third term in (43) as follows

$$\|(1-s(k))^T x(k)\| = \|(1-s(k))^T (I - \frac{1}{N}11^T)x(k)\|$$
$$\leq \sqrt{N^2 - N}\|\breve{x}(k)\|$$
$$\leq N\|\breve{x}(k)\|$$
(44)

Substituting (44) into (43) yields the desired results. The proof is thus established.

**Proof of Lemma 3.7:**

(i) $\|h\|^{\lambda,K} \leq \|\breve{h}\|^{\lambda,K} + \|h\|_J^{\lambda,K}$

Since $\|h(k)\| = \|(I - \frac{1}{N}11^T + \frac{1}{N}11^T)h(k)\|$, it thus follows that

$$\|h(k)\| \leq \|(I - \frac{1}{N}11^T)h(k)\| + \|\frac{1}{N}11^T h(k)\|$$
$$\leq \|\breve{h}(k)\| + \|h(k)\|_J$$
(45)

Multiplying the above relation with $\lambda^{-k}, k = 0, 1, \ldots$, it yields that

$$\lambda^{-k}\|h(k)\| \leq \lambda^{-k}\|\breve{h}(k)\| + \lambda^{-k}\|h(k)\|_J$$
(46)

Taking the maximum over $k = 0, \ldots, K$ on both side of (46), the desired result follows immediately.

(ii) $\|\breve{h}\|^{\lambda,K} \leq Q_1\|S^{-1}\|_{\max}\|A\|_{\max}\frac{\lambda(1-\lambda^B)}{(\lambda^B-\delta)(1-\lambda)}\|z\|^{\lambda,K} + \frac{\lambda^B}{(\lambda^B-\delta)}\sum_{t=1}^{B}\lambda^{-(t-1)}\|\breve{h}(t-1)\|$

Using $z(k) = \nabla F(x(k)) - \nabla F(x(k-1))$, the relation in Push-DIGing algorithm (4) is equivalent to

$$h(k+1) = R(k)h(k) + (S(k+1))^{-1}A(k)z(k+1)$$
(47)

Then, using Lemma 3.1, for all $k \geq B-1$, we can achieve that

$$\|\breve{h}(k+1)\| = \|\breve{J}h(k+1)\|$$
$$\leq \|\breve{J}R_B(k)h(k+1-B)\| + \|\breve{J}R_{B-1}(k)(S(k+2-B))^{-1}A(k+1-B)z(k+2-B)\| + \ldots +$$
$$+ \|\breve{J}R_1(k)(S(k))^{-1}A(k-1)z(k)\| + \|\breve{J}R_0(k)(S(k+1))^{-1}A(k)z(k+1)\|$$
$$\leq \delta\|\breve{J}h(k+1-B)\| + Q_1\|S^{-1}\|_{\max}\|A(k+1-B)z(k+2-B)\| + \ldots + Q_1\|S^{-1}\|_{\max}\|A(k-1)z(k)\|$$
$$+ Q_1\|S^{-1}\|_{\max}\|A(k)z(k+1)\|$$
$$\leq \delta\|\breve{J}h(k+1-B)\| + Q_1\|S^{-1}\|_{\max}\|A\|_{\max}\sum_{t=1}^{B}\|z(k+2-t)\|$$

(48)

where $Q_1$ is the constant defined in Lemma 3.1. Multiplying the above relation with $\lambda^{-k}, k = B-1, B, \ldots$, we have

$$\lambda^{-(k+1)} \| \breve{h}(k+1) \| \le \frac{\delta}{\lambda^B} \lambda^{-(k+1-B)} \| \breve{h}(k+1-B) \| + Q_1 \| S^{-1} \|_{\max} \| A \|_{\max} \sum_{t=1}^{B} \frac{1}{\lambda^{t-1}} \lambda^{-(k+2-t)} \| z(k+2-t) \| \quad (49)$$

In order to use the norm defined in (5), we need to take $\max_{k=0,1,\ldots,K}$, which in turn requires a relation $\lambda^{-(k+1)} \| \breve{h}(k+1) \|$ with $k = -1, 1, \ldots, B-2$. To get such a relation, without loss of generality, we assume the initial condition of (49) is

$$\lambda^{-(k+1)} \| \breve{h}(k+1) \| \le \lambda^{-(k+1)} \| \breve{h}(k+1) \|, k = -1, 1, \ldots, B \quad (50)$$

Taking the maximum over $k = -1, \ldots, B-2$ on the both sides of (49) and the maximum over $k = B-1, \ldots, K$ on the both sides of (50), by combining the obtained two relations, we immediately obtain

$$\begin{aligned}
\| \breve{h} \|^{\lambda,K} &\le \frac{\delta}{\lambda^B} \| \breve{h} \|^{\lambda,K-B} + Q_1 \| S^{-1} \|_{\max} \| A \|_{\max} \sum_{t=1}^{B} \frac{1}{\lambda^{t-1}} \| z \|^{\lambda,K+1-t} + \sum_{t=1}^{B} \lambda^{-(t-1)} \| \breve{h}(t-1) \| \\
&\le \frac{\delta}{\lambda^B} \| \breve{h} \|^{\lambda,K} + Q_1 \| S^{-1} \|_{\max} \| A \|_{\max} \| z \|^{\lambda,K} \sum_{t=1}^{B} \frac{1}{\lambda^{t-1}} + \sum_{t=1}^{B} \lambda^{-(t-1)} \| \breve{h}(t-1) \|
\end{aligned} \quad (51)$$

Therefore, we finally get

$$\| \breve{h} \|^{\lambda,K} \le Q_1 \| S^{-1} \|_{\max} \| A \|_{\max} \frac{\lambda(1-\lambda^B)}{(\lambda^B - \delta)(1-\lambda)} \| z \|^{\lambda,K} + \frac{\lambda^B}{(\lambda^B - \delta)} \sum_{t=1}^{B} \lambda^{-(t-1)} \| \breve{h}(t-1) \| \quad (52)$$

The desired result follows immediately.

(iii) $\| h \|_J^{\lambda,K} \le \frac{\| JR \|_{\max}}{\lambda} \| h \|^{\lambda,K} + \| S^{-1} \|_{\max} \| A \|_{\max} \| z \|^{\lambda,K}$

First, multiplying the relation (47) with $J = \frac{1}{N} \mathbf{1}\mathbf{1}^T$ and taking Euclidean norm at the both sides of the obtained equality and using the corresponding proprieties of Euclidean norm, we therefore have

$$\begin{aligned}
\| Jh(k+1) \| &= \| JR(k)h(k) + J(S(k+1))^{-1} A(k) z(k+1) \| \\
&\le \| JR(k)h(k) \| + \| J(S(k+1))^{-1} A(k) z(k+1) \| \\
&\le \| JR \|_{\max} \| h(k) \| + \| S^{-1} \|_{\max} \| A \|_{\max} \| z(k+1) \|
\end{aligned} \quad (53)$$

where we have employed $\| J \| = \| \frac{1}{N} \mathbf{1}\mathbf{1}^T \| = 1$ to obtain the last inequality. Multiplying the above relation with $\lambda^{-(k+1)}, k = 0, 1, \ldots$, we immediately have

$$\lambda^{-(k+1)} \| h(k+1) \|_J \le \frac{\| JR \|_{\max}}{\lambda} \lambda^{-k} \| h(k) \| + \| S^{-1} \|_{\max} \| A \|_{\max} \lambda^{-(k+1)} \| z(k+1) \| \quad (54)$$

Taking $\max_{k=0,1,\ldots,K-1}\{\cdot\}$ on the both sides of (54) gives

$$\| h \|_J^{\lambda,K} \le \frac{\| JR \|_{\max}}{\lambda} \| h \|^{\lambda,K} + \| S^{-1} \|_{\max} \| A \|_{\max} \| z \|^{\lambda,K} \quad (55)$$

From the conditions (i), (ii) and (iii), we further to show that the rest proof of Lemma 3.7. By combining the preceding items (i), (ii) and (iii), we have for all $K = 0, 1, \ldots,$ and $\lambda < 1$

$$\|h\|^{\lambda,K} \leq \frac{\|JR\|_{\max}}{\lambda}\|h\|^{\lambda,K} + \|S^{-1}\|_{\max}\|A\|_{\max}\left(1+Q_1\frac{\lambda(1-\lambda^B)}{(\lambda^B-\delta)(1-\lambda)}\right)\|z\|^{\lambda,K}$$
$$+\frac{\lambda^B}{(\lambda^B-\delta)}\sum_{t=1}^{B}\lambda^{-(t-1)}\|\breve{h}(t-1)\|$$
(56)

Rearranging the above formula and recalling $\|JR\|_{\max} < \lambda < 1$, we finally have

$$\|h\|^{\lambda,K} \leq \frac{\|S^{-1}\|_{\max}\|A\|_{\max}\left(1+Q_1\frac{\lambda(1-\lambda^B)}{(\lambda^B-\delta)(1-\lambda)}\right)}{1-\frac{\|JR\|_{\max}}{\lambda}}\|z\|^{\lambda,K} + \frac{\frac{\lambda^B}{(\lambda^B-\delta)}\sum_{t=1}^{B}\lambda^{-(t-1)}\|\breve{h}(t-1)\|}{1-\frac{\|JR\|_{\max}}{\lambda}}$$

This completes the proof.